
\baselineskip=14pt
\parskip=10pt

\magnification=\magstephalf

\def\1{{\overline{1}}}
\def\2{{\overline{2}}}
\parindent=0pt
\overfullrule=0in

\def\frac#1#2{{#1 \over #2}}

\centerline
{
 \bf How Many Rounds Should You Expect in Urn Solitaire? 
}
\rm
\bigskip
\centerline
{\it By Shalosh B. EKHAD and Doron ZEILBERGER}

{\bf Abstract}: A certain sampling process, concerning an urn with balls of two colors,
proposed in 1965 by B.E. Oakley and R.L. Perry, and discussed 
by Peter Winkler and Martin Gardner, that has an extremely simple answer for the {\it probability},
namely the constant function $\frac{1}{2}$, has a far more complicated  {\it expected duration}, that
we discover and sketch the proof of. So unlike, for example, the classical gambler's ruin problem, for which
both `probability of winning' and `expected duration' have very simple expressions, in this case the
expected number of rounds is extremely complicated, and beyond the scope of humans.

In Peter Winkler's wonderful book ([W], p. 76), the following intriguing problem, originally due to Oakely and Perry ([OP]) is posed

{\it 
\quad  ``Before you is an urn containing some green balls and some red ones (at least one of each). In Round 1
of this game, you draw a ball blindly and note its color. You then continue to draw balls (always randomly) until you get
one of the {\bf other} color; that one is then returned to the urn.

\quad  Round 2 and successive rounds are repetitions of Round 1. You play until the urn is empty; if the last ball drawn is
green, you win.

\quad How many green balls and how many red balls should you start with in the urn to maximize the probability of winning?''
}

The surprising answer is that it does not matter! The probability of winning with an urn containing $m$ green balls and $n$ red balls,
let's call it, $P(m,n)$ is {\bf always} $\frac{1}{2}$, provided that both $m$ and $n$ are strictly positive. Of course
$P(m,0)=1$ and $P(0,n)=0$.

Peter Winkler presents a really slick proof, that was suggested by Sergiu Hart, that is a bit too clever for our taste.
Let's first present another proof, less clever, yet much shorter than the original proof in [OP].

By conditioning on the number of balls, $k$, in the first round, and their color, it is readily seen that 
$P(m,n)$ satisfies the `dynamical programming' recurrence
$$
P(m,n)=\sum_{k=1}^{m} \frac{{{m} \choose {k}}}{{{m+n} \choose {k}}} \cdot \frac{n}{m+n-k} \cdot P(m-k,n) \, + \,
\sum_{k=1}^{n} \frac{{{n} \choose {k}}}{{{m+n} \choose {k}}} \cdot \frac{m}{m+n-k} \cdot P(m,n-k) \quad,
\eqno(Precurrence)
$$
subject to the boundary conditions $P(m,0)=1$, $P(0,n)=0$.

Suppose that you have no clue about the answer, can you guess it?
One way is to do {\it simulations}, and see that in about one half of the times you win.
But a better way is to write a very short procedure, that would tell you that $P(m,n)=\frac{1}{2}$ for $0 < m,n \leq 50$,
so a reasonable {\bf conjecture} is that it is $P(m,n)=\frac{1}{2}$ for {\bf all}  $m>0$ and $n>0$. In order to {\bf prove it}, all you need is
verify the {\bf trivial} binomial coefficient identity obtained by replacing $P(m ,n)$ by $\frac{1}{2}$ (when $m>0$ and $n>0$),
$P(m,0)$ by $1$, and $P(0,n)$, in other words, we have to prove, for $m,n>0$,
$$
\frac{1}{2}=\sum_{k=0}^{m-1} \frac{{{m} \choose {k-1}}}{{{m+n} \choose {k}}} \cdot \frac{n}{m+n-k} \cdot \frac{1}{2}
 \, + \,
\sum_{k=1}^{n-1} \frac{{{n} \choose {k}}}{{{m+n} \choose {k}}} \cdot \frac{m}{m+n-k} \cdot \frac{1}{2} \, + \, \frac{1}{{{m+n} \choose {n}}} \quad .
\eqno(TrivialIdentity)
$$

The Maple code that does it is

{\tt evalb(normal(convert((1/2)*sum(binomial(m,k)/binomial(m+n,k)*n/(m+n-k),k=1..m-1)+ \hfill\break
(1/2)*sum(binomial(n,k)/binomial(m+n,k)*a/(m+n-k),k=1..n-1) + 1/binomial(m+n,m),factorial)-1/2)=0); } \quad ,

that returns {\tt true}. However, one does not need Maple to verify this trivial identity, since the two sums are
{\it telescoping}, i.e.{\it gosperable}. QED!

While this proof is a bit longer than the one in [W] it is much less painful, since it involves far less
thinking!

But how many rounds should you expect Urn Solitaire to last? If you are really lucky, you can finish in one round,
and if you are really unlucky, it may take up to $m+n-1$ rounds, i.e. every ball happens to be of a different color than the
color of the previous ball, so each round is only one-ball-long. (We count rounds until it is clear who is the winner, i.e. the
game ends when there only remain balls of one color.)

The analogous {\it Dynamical Programming} recurrence, for the {\bf expected number of rounds}
(until the end of the game, i.e. when you are left with only one color), let's call it $E(m,n)$,
is {\bf very} similar, it is, for $m>0$ and $n>0$,
$$
E(m,n)= 1 \,+ \,
\sum_{k=1}^{m} \frac{{{m} \choose {k}}}{{{m+n} \choose {k}}} \cdot \frac{n}{m+n-k} \cdot E(m-k,n) \, + \,
\sum_{k=1}^{n} \frac{{{n} \choose {k}}}{{{m+n} \choose {k}}} \cdot \frac{m}{m+n-k} \cdot E(m,n-k) \quad,
\eqno(Erecurrence)
$$
subject to the boundary conditions $E(m,0)=0$ and $E(0,n)=0$.

Note that the {\bf only} difference is the extra $1$ on the right hand side. This is  reminiscent of the
classical {\bf gambler's ruin} problem in a fair casino. The probability of winning, if 
you lose a dollar with probability $\frac{1}{2}$ and win a dollar with probability $\frac{1}{2}$, and currently have $x$ dollars
and have to leave as soon as you get $0$ dollars or $N$ dollars is famously $x/N$, since it satisfies
the recurrence $p(x)\,=\,\frac{1}{2}\,(p(x-1)+p(x+1))$ subject to the boundary conditions $p(0)=0$, $p(N)=1$,
and the expected duration of the game, $E(x)$, is $x(N-x)$, since it satisfies the recurrence
$E(x)\, = \, 1\,+ \, \frac{1}{2}\,(E(x-1)+E(x+1))$ subject to the boundary conditions $E(0)=0$ and $E(N)=0$.
In this case both expressions are very simple, and can be easily guessed by humans, and then verified by humans.

Going back to our case, $P(m,n)$ is even simpler, it is {\bf identically} $\frac{1}{2}$, (for $m>0$ and $n>0$), but
alas, we are sure that no human can guess a `nice' expression for $E(m,n)$. But one of us (SBE) easily found
the next-best-thing to a closed-form formula, linear recurrences, in both $m$ and $n$ satisfied by $E(m,n)$ given by
the next theorem (and its obvious analog obtained by replacing $m$ by $n$ and using the symmetry of $E(m,n)$).

{\bf Theorem 1}.
$$
2\, ( m+2 )  ( m+1 )  ( 2\,mn+2\,m+7\,n+9 )  ( 3+m) E ( m,n)  
$$
$$
- ( 3+m )  ( m+2)  ( 14\,{m}^{2}n+6\,m{n}^{2}+14\,{m}^{2}+83\,mn+21\,{n}^{2}+91\,m+123\,n+128 ) E( m+1,n)
$$
$$ 
+ (3+m) ( 18\,{m}^{3}n+16\,{m}^{2}{n}^{2}+2\,m{n}^{3}+18\,{m}^{3}+167\,{m}^{2}n+94\,m{n}^{2}
$$
$$
+7\,{n}^{3}+169\,{m}^{2}+521\,mn+132\,{n}^{2}+511\,m+537\,n+508) E( m+2,n) 
$$
$$
+ ( -10\,{m}^{4}n-14\,{m}^{3}{n}^{2}-4\,{m}^{2}{n}^{3}-10\,{m}^{4}-135\,{m}^{3}n-129\,{m}^{2}{n}^{2}-24\,m{n}^{3}
-131\,{m}^{3}-
$$
$$
684\,{m}^{2}n-393\,m{n}^{2}-34\,{n}^{3}-639\,{m}^{2}-1535\,mn-396\,{n}^{2}-1380\,m-1286\,n-1116) E ( m+3,n) 
$$
$$
+ ( 3+m )  ( 2\,mn+2\,m+5\,n+7 )  ( n+4+m )  ( n+3+m ) E ( m+4,n ) =0 \quad .
\eqno(NiceRecurrence)
$$

{\bf Sketch of Proof}: The proof is similar to the one for $P(m,n)=\frac{1}{2}$, `plug-in' the recurrence
$(NiceRecurrence)$ into $(Erecurrence)$.
Now one gets an identity in the {\bf holonomic ansatz}, that could be rigorously proved using Christoph  Koutschan's
amazing Mathematica package ([K1][K2]). But since a proof {\it exists}, we can get a semi-rigorous proof just by checking it
for sufficiently many values of $E(m,n)$ and that is what we did.

Note that using $(NiceRecurrence)$ (and the analogous one w.r.t. $n$) one can compute $E(m,n)$ in {\it linear time} and {\it constant memory},
while using the $(Erecurrence)$ takes both quadratic time and memory.
Using the holonomic ansatz, we  derived a third-order linear recurrence with polynomial coefficients for the diagonal
sequence $E(n,n)$.

{\bf Theorem 2}: The expected number of rounds in Urn Solitaire with $n$ balls of each color satisfies the recurrence
$$
-2\, \left( 18\,{n}^{4}+159\,{n}^{3}+528\,{n}^{2}+779\,n+428 \right)  \left( n+1 \right) ^{2} \left( n+2 \right) E \left( n,n \right) 
$$
$$
+ \left( n+2 \right)  \left( 216\,{n}^{6}+2358\,{n}^{5}+10485\,{n}^{4}+24174\,{n}^{3}+30251\,{n}^{2}+19276\,n+4800 \right) E \left( n+1,n+1 \right) 
$$
$$
+ \left( -324\,{n}^{7}-4032\,{n}^{6}-21015\,{n}^{5}-59334\,{n}^{4}-97813\,{n}^{3}-93898\,{n}^{2}-48288\,n-10080 \right) E \left( n+2,n+2 \right) 
$$
$$
+2\, \left( n+2 \right)  \left( 18\,{n}^{4}+87\,{n}^{3}+159\,{n}^{2}+128\,n+36 \right)  \left( 2\,n+5 \right) ^{2}E \left( n+3 , n+3 \right) =0 \quad ,
$$
subject to the initial conditions $E(1,1)=1$, $E(2,2)=\frac{17}{9}$, $E(3,3)=\frac{143}{50}$.

{\bf  Maple Packages}

Everything in this paper was found using the Maple packages {\tt UrnSolitaire.txt} and {\tt GuessHolo2.txt}, available, along with sample input 
and output files, from the webpage of this article

{ \tt http://sites.math.rutgers.edu/\~{}zeilberg/mamarim/mamarimhtml/urn.html} \quad .

{\bf Simple Urn Solitaire}

Suppose that you do {\bf not} return to the urn the balls that end a round, but otherwise define a round the same way, then
instead of $P(m,n)=\frac{1}{2}$ the probability is, of course $\frac{m}{m+n}$, and the expected number of
rounds, let's call it $F(m,n)$ satisfies the dynamical programming recurrence 
for $m>0$ and $n>0$,
$$
F(m,n)= 1 \,+ \,
\sum_{k=1}^{m} \frac{{{m} \choose {k}}}{{{m+n} \choose {k}}} \cdot \frac{n}{m+n-k} \cdot F(m-k,n-1) \, + \,
\sum_{k=1}^{n} \frac{{{n} \choose {k}}}{{{m+n} \choose {k}}} \cdot \frac{m}{m+n-k} \cdot F(m-1,n-k) \quad,
\eqno(Frecurrence)
$$
subject to the boundary conditions $F(m,0)=0$ and $F(0,n)=0$. For the recurrence for the expectation
(and also for the variance!) for the diagonal sequence  $F(n,n)$ (and the corresponding one for the variance)
see the output file

{ \tt http://sites.math.rutgers.edu/\~{}zeilberg/tokhniot/oUrnSolitaire2.txt} \quad .

{\bf Conclusion} 

Peter Winkler's book is full of challenging problems that (smart) humans can do all
by themselves, but take any of these problems, and tweak it ever-so-slightly, and then humans are
hopeless, but luckily they can ask computer-kind to do their work.

{\bf References}

[G] Martin Gardner, ``{\it The Colossal Book of Short Puzzles and Problems}'', W.W. Norton \& Co, 2005.

[K1] Christoph Koutschan,
{\it ``Advanced Applications of the Holonomic Systems Approach''},  RISC-Linz, Johannes Kepler University. PhD Thesis. September 2009.

[K2] Christoph Koutschan, {\it HolonomicFunctions (User's Guide). Technical report no. 10-01 in RISC Report Series, 
Research Institute for Symbolic Computation (RISC), Johannes Kepler University Linz, Austria}.  January 2010.

[OP] B.E. Oakley and R.L. Perry, {\it A sampling process}, The Mathematical Gazette {\bf 49} No. 367 (Feb. 1965), 42-44.

[W] Peter Winkler, ``{\it Mathematical Mind-Benders},'' A.K. Peters/CRC Press, 2007. 

\vfill\eject

\bigskip
\bigskip
\hrule
\bigskip
Doron Zeilberger, Department of Mathematics, Rutgers University (New Brunswick), Hill Center-Busch Campus, 110 Frelinghuysen
Rd., Piscataway, NJ 08854-8019, USA. \hfill \break
DoronZeil at gmail dot com  \quad ;  \quad {\tt http://sites.math.rutgers.edu/\~{}zeilberg/} \quad .
\bigskip
\hrule
\bigskip
Shalosh B. Ekhad, c/o D. Zeilberger, Department of Mathematics, Rutgers University (New Brunswick), Hill Center-Busch Campus, 110 Frelinghuysen
Rd., Piscataway, NJ 08854-8019, USA.
\bigskip
\hrule

\bigskip
Exclusively published in The Personal Journal of Shalosh B. Ekhad and Doron Zeilberger  \hfill \break
{ \tt http://sites.math.rutgers.edu/\~{}zeilberg/pj.html} \quad  and {\tt arxiv.org} \quad . 
\bigskip
\hrule
\bigskip

{\bf  Written:  Jan. 4, 2018.} 

\end